\documentclass[12pt]{article}
\usepackage{amsmath}
\usepackage{color}
\usepackage{fancyhdr}
\usepackage{verbatim}
\usepackage{amsfonts,longtable,mathrsfs,amsfonts,longtable,mathrsfs,longtable,mathrsfs,mathtools, float, amsmath,epsfig}
\usepackage[numbers,sort&compress]{natbib}
\usepackage{hyperref}
\usepackage{sectsty} 

\newtheorem{Theorem}{Theorem}

\newtheorem{Proof}{Proof}

\def\R{\hbox{{\rm I}\kern-0.2em{\rm R}\kern0.2em}}

\def\bn{\begin{equation}}
\def\en{\end{equation}}
\def\bny{\begin{eqnarray}}
\def\eny{\end{eqnarray}}
\def\be{\begin{eqnarray*}}
\def\ee{\end{eqnarray*}}
\def\bc{\begin{center}}
\def\ec{\end{center}}

\def\({\left(}
\def\){\right  )}
\def\[{\left[}
\def\]{\right]}
\def\bc{\begin{center}}
\def\ec{\end{center}}

\newtheorem{dfn}{Definition}[section]
\newtheorem{thm}{Theorem}[section]
\newtheorem{rem}{Remark}[section]
\newtheorem{pro}{Proposition}[section]

\newtheorem{cor}{Corollary}[section]
\newtheorem{lem}{Lemma}[section]
\newtheorem{exm}{Example}[section]

\def\bn{\begin{equation}}
\def\en{\end{equation}}
\def\bny{\begin{eqnarray}}
\def\eny{\end{eqnarray}}
\def\be{\begin{eqnarray*}}
\def\ee{\end{eqnarray*}}
\def\bdn{\begin{dfn}}
\def\edn{\end{dfn}}
\def\btm{\begin{thm}}
\def\etm{\end{thm}}
\def\bpf{\begin{proof}}
\def\epf{\end{proof}}
\def\bpn{\begin{pro}}
\def\epn{\end{pro}}
\def\brk{\begin{rem}}
\def\erk{\end{rem}}
\def\bcy{\begin{cor}}
\def\ecy{\end{cor}}
\def\blm{\begin{lem}}\def\elm{\end{lem}}
\def\bex{\begin{exm}}
\def\eex{\end{exm}}

 \def\R{{\hat R}}
\begin{document}

\bc {Method of Lie Symmetry for analytical solutions, periodicity and attractivity of a family of tenth-order difference equations}\ec
\bc
\bc{M. Folly-Gbetoula\footnote{Corresponding author (Mensah.Folly-Gbetoula@wits.ac.za)} and Kwassi Anani\\
}
\ec
School of Mathematics, University of the Witwatersrand, 2050, Johannesburg, South Africa.\\
Département de Mathématiques, Faculté des Sciences, Université de Lomé, 1515, Lomé, Togo.
\ec
\begin{abstract}
\noindent Symmetry is a powerful tool for finding analytical solutions to differential equations, both partial and ordinary, via the similarity variables or via the invariance of the equation under group transformations. It is the largest group of transformations that leaves the differential  equation invariant. It is now known that this differential equation method plays the same role when it comes to the study of difference equations.    Difference equations can be used to model various phenomena where the changes occur in discrete manner. The use of symmetries on recurrence equations, usually, leads to reductions of order and hence ease the process of finding their solutions. One of the aims of this work is to employ symmetries to generalize some results in the literature. We present new generalized formula solutions of a class of difference equations and we investigate the periodicity and behavior of theses solutions.
\end{abstract}
\section{Introduction} \setcounter{equation}{0}
Difference equations have attracted a lot of attention. In practice, they are are used to model phenomena in which the variable is discrete. The idea of symmetries was poineered by Sophus Lie (1842-1899), where he applied it to differential equations.  
In 1987, Maeda demonstrated that some modified version of Lie's approach could also be applied to  find solutions of ordinary difference equations \cite{Maeda}.\\
\noindent In this study, we use Lie symmetry analysis to find solutions of a certain class of difference equations. Our work is inspired by the work of A. Sanbo and Elsayed \cite{SE}, where the authors studied the difference equations 
\begin{equation}\label{un}
x_{n+1}=\frac{x_{n-9}}{\pm1\pm x_{n-1}x_{n-3}x_{n-5}x_{n-7}x_{n-9}}, \quad n=0, 1, \cdots.
\end{equation}
In their work, Sanbo and Elsayed \cite{SE}, obtained expressions for solutions of \eqref{un} and analyzed the dynamical behavior of these solutions. The solutions were proved by mathematical induction.\\
Our approach is different from the method Sanbo and Elsayed used to find  the solutions. We use Lie symmetry analysis to solve the generalized difference equation: 
\begin{equation}\label{un2}
    x_{n+1}=\frac{x_{n-9}}{a_n+b_nx_{n-1}x_{n-3}x_{n-5}x_{n-7}x_{n-9}},
\end{equation}
where $a_n$ and $b_n$ are real numbers.
Clearly, this is more general than what was considered by Sanbo and Elsayed. For similar work on difference equations, see \cite{EM, El1,tf}.
\subsection{Symmetries for difference equations}
\noindent Consider the tenth-order ordinary recurrence equation
\begin{equation}\label{5}
    x_{n+10}=\Omega(n,x_n,x_{n+2},x_{n+4},x_{n+6},x_{n+8}), \quad \frac{\partial\Omega}{\partial x_n}\neq 0,
\end{equation}
where $\Omega$ is some smooth function and $n$ an independent variable. The solution of \eqref{5} may be expressed as 
\begin{equation}
    x_n=f(n,c_1,\cdots,c_{10}).
\end{equation}
\\
\begin{dfn}
	The forward shift operator is given by
	\begin{equation}
	S:n\mapsto n+1, \qquad S^{i}x_n= x_{n+i}.
	\end{equation}
\end{dfn}
    The tenth-order ordinary recurrence equation \eqref{5} admits a symmetry generator $X$ given as
    \begin{align}\label{Ngener}
X=&Q\frac{\partial}{ \partial x_n}+ Q(n+2, x_{n+2})\frac{\partial}{ \partial x_{n+2}}+ Q(n+4, x_{n+4})\frac{\partial}{ \partial x_{n+4}}+\nonumber \\& Q(n+6, x_{n+6})\frac{\partial}{ \partial x_{n+6}}+ Q(n+8, x_{n+8})\frac{\partial}{ \partial x_{n+8}}
\end{align}
that satisfies the symmetry condition
\begin{equation}\label{LSC}
\mathcal{S}^{(10)} Q- X \Omega=0.
\end{equation}
The function $Q=Q(n, x_n)$ is called the characteristic of the group of transformations. For more details on this, see \cite{Hydon,M,mmd,FK,PH}.
\section{Main results}
We will use Lie point symmetry in this section to obtain solutions to the tenth-order ordinary recurrence equation \eqref{un2}. We consider the  equivalent tenth-order ordinary recurrence equation
\begin{equation}\label{shift}
    x_{n+10} = \frac{x_n}{A_n+B_nx_nx_{n+2}x_{n+4}x_{n+6}x_{n+8}}
\end{equation}
of \eqref{un2}, where $A_n=a_{n+9}$ and $B_n=b_{n+9}$. Applying the condition \eqref{LSC} to \eqref{shift}, we obtain
\begin{align}\label{a1}
&Q(n+10, \Omega)-\bigg[Q(n, x_n)\frac{\partial \Omega}{\partial x_n}+Q(n+2, x_{n+2})\frac{\partial \Omega}{\partial x_{n+2}}+Q(n+4, x_{n+4})\frac{\partial \Omega}{\partial x_{n+4}}+\nonumber\\&+Q(n+6, x_{n+6})\frac{\partial \Omega}{\partial x_{n+6}}+Q(n+8, x_{n+8})\frac{\partial \Omega}{\partial x_{n+8}}\bigg]=0
\end{align}
where $\Omega$ is the right-hand side expression in \eqref{shift}). At this step, any convenient operator can be applied to \eqref{a1}, provided that it leads to a reduction in the recurrence equation.
Applying the operator  \begin{align}\label{D}
	L=&\frac{\partial}{\partial x_n}+\frac{A_n}{B_nx_n^2x_{n+4}x_{n+6}x_{n+8}}\frac{\partial\quad\quad}{\partial x_{n+2}}
\end{align}  to \eqref{a1}, we get (after simplification)
\begin{align}\label{a13}
& -\alpha_n Q'(n, x_n)+ \frac{2\alpha_n}{x_n}Q(n, x_n)+ \alpha_n Q'(n+2, x_{n+2}) +\frac{\alpha_n}{x_{n+4}}Q(n+4, x_{n+4})\nonumber \\& +\frac{\alpha_n}{x_{n+6}}Q(n+6, x_{n+6}) +\frac{\alpha_n}{x_{n+8}}Q(n+8, x_{n+8})=0,
\end{align}
where $Q'(n,x_n)$ denotes the derivative of $Q$ with respect to $x_n$. To overcome the challenge of dealing with different arguments, we differentiate \eqref{a13} with respect to $x_n$ and obtain the following:
\begin{eqnarray}\label{f1}
    &\dfrac{d\quad}{dx_n}\left(-Q'(n, x_n)+\frac{2}{x_n}Q(n, x_n)\right)= 0.
\end{eqnarray}
Thus the solution is \begin{equation}
Q(n,x_n)=\gamma_n x_n^2+\alpha_nx_n+\beta_n,
\end{equation}
where $\alpha_n$, $\beta_n$ and $\gamma_n$ are some functions of $n$.\\
By substituting $Q(n,x_n)$ in the symmetry condition \eqref{LSC}, we obtain that $\gamma_n=\beta_n=0$, and $\alpha_n$ must satisfy
 \begin{equation}\label{f4}
 \alpha_{n}{+}\alpha_{n+2}{+}\alpha_{n+4}{+}\alpha_{n+6}{+}\alpha_{n+8}=0.
 \end{equation}
Thus the symmetry generator takes the form  
 \begin{equation}\label{f3}
    X=Q(n, x_n)\frac{\partial}{\partial x_n}= \alpha_nx_n\frac{\partial}{\partial x_n},
\end{equation}
where $\alpha_n$ satisfies \eqref{f4}.
The solutions of \eqref{f4} are 
\begin{equation}\label{alpha}
\alpha_n=e^{\frac{i(2kn\pi)}{10}},
\end{equation}
$k=1,2,3,4,6,7,8,9$.
It follows from \eqref{f3} and \eqref{alpha} that 
\begin{equation}\label{f3'}
X_k= e^{i\frac{2kn\pi}{10}}x_n\frac{\partial}{\partial x_n},
\end{equation}
$k=1,2,3,4,6,7,8,9$, are symmetries of \eqref{shift}. We introduce the canonical coordinate defined as
\begin{align}
S_n =& \int_{}^{}\frac{1}{Q(n,x_n)}dx_n=\int_{}^{}\frac{1}{\alpha_nx_n}dx_n,
\end{align}
we get
\begin{equation}
    S_n\alpha_n=\ln|x_n|.
\end{equation}
Let $\Tilde{F_n}=S_n\alpha_n+S_{n+2}\alpha_{n+4}+S_{n+6}\alpha_{n+6}+S_{n+8}\alpha_{n+8}$
and 
$F_n =e^{{-{\Tilde{F_n}}}}.
$
Hence, we have $\Tilde{F_n}=\ln{(x_nx_{n+2}x_{n+4}x_{n+6}x_{n+8})}$
and
\begin{align}\label{V}
F_n =\frac{1}{x_nx_{n+2}x_{n+4}x_{n+6}x_{n+8}}.
\end{align}
Shifting equation \eqref{V} twice and substituting the expression of  $x_{n+10}$ given in \eqref{shift}, we obtain
\begin{align}\label{V2}
    F_{n+2}=A_nF_n+B_n.
\end{align}
By iterating (\ref{V2}), we are led to:
\begin{align}\label{I}
    F_{2n+j}=F_j\left(\prod_{t=0}^{n-1}A_{2t+j}\right)+\sum_{i=0}^{n-1}\left(B_{2i+j}\prod_{k_2=i+1}^{n-1}A_{2k_2+j}\right),\quad j=0,1.
\end{align}
Also,
\begin{align}
    x_{n+10}=\frac{F_n}{F_{n+2}}x_n
\end{align}
whose iteration yields
\begin{align}\label{I1}
x_{10n+k}=&x_k\Bigg({\prod_{s=0}^{n-1}\frac{F_{10s+k}}{F_{10s+k+2}}}\Bigg),\quad k=0,1,2,\cdots,9\nonumber\\
=&x_k\Bigg({\prod_{s=0}^{n-1}\frac{F_{10s+2\lfloor \frac{k}{2}\rfloor+\tau{(k)}}}{F_{10s+2\lfloor \frac{k}{2}\rfloor+\tau{(k)}+2}}}\Bigg)\nonumber\\
=&x_k\Bigg({\prod_{s=0}^{n-1}\frac{F_{2(5s+\lfloor \frac{k}{2}\rfloor)+\tau{(k)}}}{F_{2(5s+1+\lfloor \frac{k}{2}\rfloor)+\tau{(k)}}}}\Bigg)
\end{align}
since we can always write any integer in the form $i=2\lfloor \frac{i}{2}\rfloor +\tau(i)$, where $\tau(i)$ is the remainder when $a$ is divided by $2$ and $\lfloor \cdot \rfloor $ is the floor function.
\noindent Substituting equation (\ref{I}) into equation (\ref{I1}), we get
\begin{align}\label{anbn}
x_{10n+k}=&x_k{\prod_{s=0}^{n-1}\frac{F_{\tau(k)}\left(\prod\limits_{t=0}^{5s+\lfloor \frac{k}{2}\rfloor-1}A_{2t+{\tau(k)}}\right)+\sum\limits_{i=0}^{5s+\lfloor \frac{k}{2}\rfloor-1}\left(B_{2i+{\tau(k)}}\prod\limits_{k_2=i+1}^{5s+\lfloor \frac{k}{2}\rfloor-1}A_{2k_2+{\tau(k)}}\right)}{F_{\tau(k)}\left(\prod\limits_{t=0}^{5s+\lfloor \frac{k}{2}\rfloor}A_{2t+{\tau(k)}}\right)+\sum\limits_{i=0}^{5s+\lfloor \frac{k}{2}\rfloor}\left(B_{2i+{\tau(k)}}\prod\limits_{k_2=i+1}^{5s+\lfloor \frac{k}{2}\rfloor}A_{2k_2+{\tau(k)}}\right)}}\nonumber\\
=&x_k{\prod_{s=0}^{n-1}\frac{\left(\prod\limits_{t=0}^{5s+\lfloor \frac{k}{2}\rfloor-1}A_{2t+{\tau(k)}}\right)+\sum\limits_{i=0}^{5s+\lfloor \frac{k}{2}\rfloor-1}\left(\frac{B_{2i+{\tau(k)}}}{F_{\tau(k)}}\prod\limits_{k_2=i+1}^{5s+\lfloor \frac{k}{2}\rfloor-1}A_{2k_2+{\tau(k)}}\right)}{\left(\prod\limits_{t=0}^{5s+\lfloor \frac{k}{2}\rfloor}A_{2t+{\tau(k)}}\right)+\sum\limits_{i=0}^{5s+\lfloor \frac{k}{2}\rfloor}\left(\frac{B_{2i+{\tau(k)}}}{F_{\tau(k)}}\prod\limits_{k_2=i+1}^{5s+\lfloor \frac{k}{2}\rfloor}A_{2k_2+{\tau(k)}}\right)}},
\end{align}
where $1/F_{\tau(k)}=x_{\tau(k)}x_{\tau(k)+2}x_{\tau(k)+4}x_{\tau(k)+6}x_{\tau(k)+8}$.
\section{Case where $A_n = A$ and $B_n=B$}
Let $A_n=A$ and $B_n=B$ in equation (\ref{anbn}). We obtain
\begin{align}\label{I3}
x_{10n+k}=&x_k\left({\prod_{s=0}^{n-1}\frac{A^{5s+\lfloor \frac{k}{2}\rfloor}+\frac{B}{F_{\tau(k)}}\sum\limits_{s=0}^{5s+\lfloor \frac{k}{2}\rfloor-1}A^{s}}{A^{5s+1+\lfloor \frac{k}{2}\rfloor}+\frac{B}{F_{\tau(k)}}\sum\limits_{s=0}^{5s+\lfloor \frac{k}{2}\rfloor}A^{s}}}\right)\end{align}\begin{align}
=&x_k\left({\prod_{s=0}^{n-1}\frac{A^{5s+\lfloor \frac{k}{2}\rfloor}+{B}{x_{\tau(k)}x_{\tau(k)+2}x_{\tau(k)+4}x_{\tau(k)+6}x_{\tau(k)+8}}\sum\limits_{s=0}^{5s+\lfloor \frac{k}{2}\rfloor-1}A^{s}}{A^{5s+1+\lfloor \frac{k}{2}\rfloor}+{B}{x_{\tau(k)}x_{\tau(k)+2}x_{\tau(k)+4}x_{\tau(k)+6}x_{\tau(k)+8}}\sum\limits_{s=0}^{5s+\lfloor \frac{k}{2}\rfloor}A^{s}}}\right)
\end{align}\begin{align}
=& \begin{cases}
x_k{\prod\limits_{s=0}^{n-1}\frac{1+{B}(5s+\lfloor \frac{k}{2}\rfloor){x_{\tau(k)}x_{\tau(k)+2}x_{\tau(k)+4}x_{\tau(k)+6}x_{\tau(k)+8}}}{1+{B}(5s+\lfloor \frac{k}{2}\rfloor+1){x_{\tau(k)}x_{\tau(k)+2}x_{\tau(k)+4}x_{\tau(k)+6}x_{\tau(k)+8}}}}, \text{when} \; A=1;
\\ ~\\
x_k{\prod\limits_{s=0}^{n-1}\frac{A^{5s+\lfloor \frac{k}{2}\rfloor}+{\frac{B(1-A^{5s+\lfloor \frac{k}{2}\rfloor})}{1-A}}{x_{\tau(k)}x_{\tau(k)+2}x_{\tau(k)+4}x_{\tau(k)+6}x_{\tau(k)+8}}}{A^{5s+1+\lfloor \frac{k}{2}\rfloor}+\frac{B(1-A^{5s+1+\lfloor \frac{k}{2}\rfloor})}{1-A}{x_{\tau(k)}x_{\tau(k)+2}x_{\tau(k)+4}x_{\tau(k)+6}x_{\tau(k)+8}}}}, \text{when} \; A\neq 1;
\end{cases}\label{sol1}
\end{align}
for $k = 0,1,2\cdots,9$.
\subsection{Case where $A = -1$}
For the case $A=-1$, we have 
\begin{align}
x_{10n+k}=x_k\left({\prod\limits_{s=0}^{n-1}\frac{{(-1)}^{s}+{\frac{B((-1)^{\lfloor \frac{k}{2}\rfloor}-(-1)^{s})}{2}}{x_{\tau(k)}x_{\tau(k)+2}x_{\tau(k)+4}x_{\tau(k)+6}x_{\tau(k)+8}}}{-(-1)^{s}+\frac{B((-1)^{\lfloor \frac{k}{2}\rfloor}+(-1)^{s})}{2}{x_{\tau(k)}x_{\tau(k)+2}x_{\tau(k)+4}x_{\tau(k)+6}x_{\tau(k)+8}}}}\right).
\end{align}
The above implies that 
\begin{subequations}\label{a=-1p}:
	\begin{itemize}
\item For $k=0,1$;
\begin{align}
x_{10n+k}=\begin{cases}
x_k   \hspace{3.9cm}\text{ if } \quad n\quad  \text{even} \\
\frac{x_k}{-1+Bx_{k}x_{k+2}x_{k+4}x_{k+6}x_{k+8}}\quad  \text{if}\quad   n\quad  \text{odd}
\end{cases},
\end{align}
\item For $k=2,3$;
\begin{align}
x_{10n+k}=\begin{cases}
x_k \quad\; \hspace{5cm}\text{ if } \quad n\quad  \text{even} \\
{x_k}({-1+Bx_{k-2}x_{k}x_{k+2}x_{k+4}x_{k+6}})\quad  \text{if} \quad n\quad \text{odd}
\end{cases},
\end{align}
\item For $k=4,5$;
\begin{align}
x_{10n+k}=\begin{cases}
x_k   \hspace{3.8cm} \text{ if } \quad n\quad  \text{even} \\
\frac{x_k}{-1+Bx_{k-4}x_{k-2}x_{k}x_{k+2}x_{k+4}}\quad  \text{if} \quad n\quad \text{odd}
\end{cases},
\end{align}
\item For $6,7$;
\begin{align}
x_{10n+k}=\begin{cases}
x_k  \quad\, \hspace{5cm}\; \text{ if } \quad n\quad  \text{even} \\
{x_k}{(-1+Bx_{k-6}x_{k-4}x_{k-2}x_{k}x_{k+2})}\quad  \text{if} \quad n\quad \text{odd}
\end{cases},
\end{align}
\item For $8,9$;
\begin{align}
x_{10n+k}=\begin{cases}
x_k \quad  \quad \hspace{5cm}\text{ if } \quad n\quad  \text{even} \\
\frac{x_k}{(-1+Bx_{k-8}x_{k-6}x_{k-4}x_{k-2}x_{k})}\quad  \text{if} \quad n\quad \text{odd}
\end{cases}.
\end{align}
\end{itemize}
\end{subequations}
\subsection{Known cases in the literature}
Recall that we shifted equation \eqref{un2} forward $9$ times to obtain \eqref{shift}, the solution of which is provided by \eqref{sol1}. We reverse the equations in \eqref{sol1} $9$ times to obtain the solution of the difference equation \eqref{un2}, which is given as

\begin{align}\label{shift1}
x_{10n+k-9}=& 
x_{k-9}{\prod\limits_{s=0}^{n-1}\frac{1+{B}(5s+\lfloor \frac{k}{2}\rfloor){x_{\tau(k)-9}x_{\tau(k)-7}x_{\tau(k)-5}x_{\tau(k)-3}x_{\tau(k)-1}}}{1+{B}(5s+\lfloor \frac{k}{2}\rfloor+1){x_{\tau(k)-9}x_{\tau(k)-7}x_{\tau(k)-5}x_{\tau(k)-3}x_{\tau(k)-1}}}},
\end{align}
 if $A = 1$; and 
\begin{align}\label{shift2}
x_{10n+k-9}=x_{k-9}{\prod\limits_{s=0}^{n-1}\frac{A^{5s+\lfloor \frac{k}{2}\rfloor}+{\frac{B(1-A^{5s+\lfloor \frac{k}{2}\rfloor})}{1-A}}{x_{\tau(k)-9}x_{\tau(k)-7}x_{\tau(k)-5}x_{\tau(k)-3}x_{\tau(k)-1}}}{A^{5s+1+\lfloor \frac{k}{2}\rfloor}+\frac{B(1-A^{5s+1+\lfloor \frac{k}{2}\rfloor})}{1-A}{x_{\tau(k)-9}x_{\tau(k)-7}x_{\tau(k)-5}x_{\tau(k)-3}x_{\tau(k)-1}}}},
\end{align}
if $A\neq 1$.\par 
\begin{subequations}\label{change}
If we set \begin{align}
j=9-k, 
\end{align}
then 
\begin{align}
\lfloor \frac{j}{2}\rfloor=4-\lfloor \frac{k}{2}\rfloor
\end{align}
and 
\begin{align}
  \tau{(j)}=1-\tau(k)
  \end{align}
 \end{subequations} 
   for $k=0,1,\dots,9$. Therefore, from \eqref{shift1}, \eqref{shift2} and \eqref{change}, we have 
\begin{align}\label{shift1'}
x_{10n-j}=& 
x_{-j}{\prod\limits_{s=0}^{n-1}\frac{1+{B}(5s+4-\lfloor \frac{j}{2}\rfloor){x_{-\tau(j)-8}x_{-\tau(j)-6}x_{-\tau(j)-4}x_{-\tau(j)-2}x_{-\tau(j)}}}{1+{B}(5s+5-\lfloor \frac{j}{2}\rfloor){x_{-\tau(j)-8}x_{-\tau(j)-6}x_{-\tau(j)-4}x_{-\tau(j)-2}x_{-\tau(j)}}}}, 
\end{align}
if $A=1$; and 
\begin{align}\label{shift2'}
x_{10n-j}=x_{-j}{\prod\limits_{s=0}^{n-1}\frac{A^{5s+4-\lfloor \frac{j}{2}\rfloor}+{\frac{B(1-A^{5s+4-\lfloor \frac{j}{2}\rfloor})}{1-A}}{x_{-\tau(j)-8}x_{-\tau(j)-6}x_{-\tau(j)-4}x_{-\tau(j)-2}x_{-\tau(j)}}}{A^{5s+5-\lfloor \frac{j}{2}\rfloor}+\frac{B(1-A^{5s+5-\lfloor \frac{j}{2}\rfloor})}{1-A}{x_{-\tau(j)-8}x_{-\tau(j)-6}x_{-\tau(j)-4}x_{-\tau(j)-2}x_{-\tau(j)}}}},
\end{align}
if $A\neq 1$.
\subsubsection{The case: $A=1$ and $B=1$}
Set $M_j=5-\lfloor\frac{j}{2}\rfloor$, $P_j=\prod\limits_{k=0}^{4}a_{\tau(j)+2k}$ and $x_{-j}=a_j$. We have
\begin{equation}
    P_j=\prod\limits_{k=0}^{4}a_{\tau(j)+2k}
        \\= x_{-\tau(j)-8}x_{-\tau(j)-6}x_{-\tau(j)-4}x_{-\tau(j)-2}x_{-\tau(j)}.
\end{equation}
Then, equation \ref{shift1'} becomes 
\begin{equation}
    x_{10n-j} = a_j\prod\limits_{s=0}^{n-1}\bigg(\frac{1+(5s+M_j-1)P_j}{1+(5s+M_j)P_j}\bigg),
\end{equation}
which is the same as Theorem 2 in \cite{SE}.
\subsubsection{Case: $A=1$ and $B=-1$}
In this case, equation \eqref{shift1'} is given as
\begin{equation}
    x_{10n-j} = a_j\prod\limits_{s=0}^{n-1}\bigg(\frac{1-(5s+M_j-1)P_j}{1-(5s+M_j)P_j}\bigg),
\end{equation}
which is the same as Theorem 5 in \cite{SE}.

\subsubsection{Case: $A=-1$ and $B=1$}

In this case, equation \eqref{shift2'} becomes 
\begin{equation}
    x_{10n-j}=a_j\prod\limits_{s=0}^{n-1}\frac{(-1)^{5s+M_j-1}+\frac{(1-(-1)^{5s+M_j-1})}{2}P_j}{(-1)^{5s+M_j}+\frac{(1-(-1)^{5s+M_j})}{2}P_j},
\end{equation}
which was obtained in Theorem 7 of \cite{SE}.
\subsubsection{Case: $A=-1$ and $B=-1$}
As above, in this case, equation \eqref{shift2'} yields
\begin{equation}
    x_{10n-j}=a_j\prod\limits_{s=0}^{n-1}\frac{(-1)^{5s+M_j-1}+\frac{(-1-(-1)^{5s+M_j-2})}{2}P_j}{(-1)^{5s+M_j}+\frac{(-1-(-1)^{5s+M_j-1})}{2}P_j},
\end{equation}
which was obtained in Theorem 10 of \cite{SE}.
\section{Periodic property and stability}
In this section, we provide sufficient conditions for existence of periodic solutions with periods $1$, $5$, $10$ and $20$. We further prove that the equilibrium points are non-hyperbolic.
\begin{Theorem}\label{0}
	Let $x_n$ be a solution of 
	{\footnotesize 
		\begin{equation}\label{Teq1cstA1}
			x_{n+10}=\frac{x_{n}}{-1+Bx_{n}x_{n+2}x_{n+4}x_{n+6}x_{n+8}}
	\end{equation}}
	with initial conditions $x_i, \; i=0,\dots, 9$; and $B$ is a non-zero constant. Then we have a $20$-periodic solution.
\end{Theorem}
\begin{Proof}
Clearly, $A=-1$. It follows from \eqref{a=-1p} that the $20$-periodic solutions are given by 
\begin{align}
&	\dots, x_0, x_1, x_2, x_3, x_4, x_5, x_6, x_7, x_8, x_9, \frac{x_0}{-1+Bx_0x_2x_4x_6x_8},\nonumber\\ &\frac{x_1}{-1+Bx_1x_3x_5x_7x_9}, {x_2}{(-1+Bx_0x_2x_4x_6x_8)}, {x_3}{(-1+Bx_1x_3x_5x_7x_9)},\nonumber\\& \frac{x_4}{-1+Bx_0x_2x_4x_6x_8}, \frac{x_5}{-1+Bx_1x_3x_5x_7x_9}, {x_6}{(-1+Bx_0x_2x_4x_6x_8)},\nonumber\\& {x_7}{(-1+Bx_1x_3x_5x_7x_9)}, \frac{x_8}{-1+Bx_0x_2x_4x_6x_8}, \frac{x_9}{-1+Bx_1x_3x_5x_7x_9}, \cdots.
\end{align}
\end{Proof}
\begin{figure}[H]
	\centering
	\includegraphics[scale=0.2]{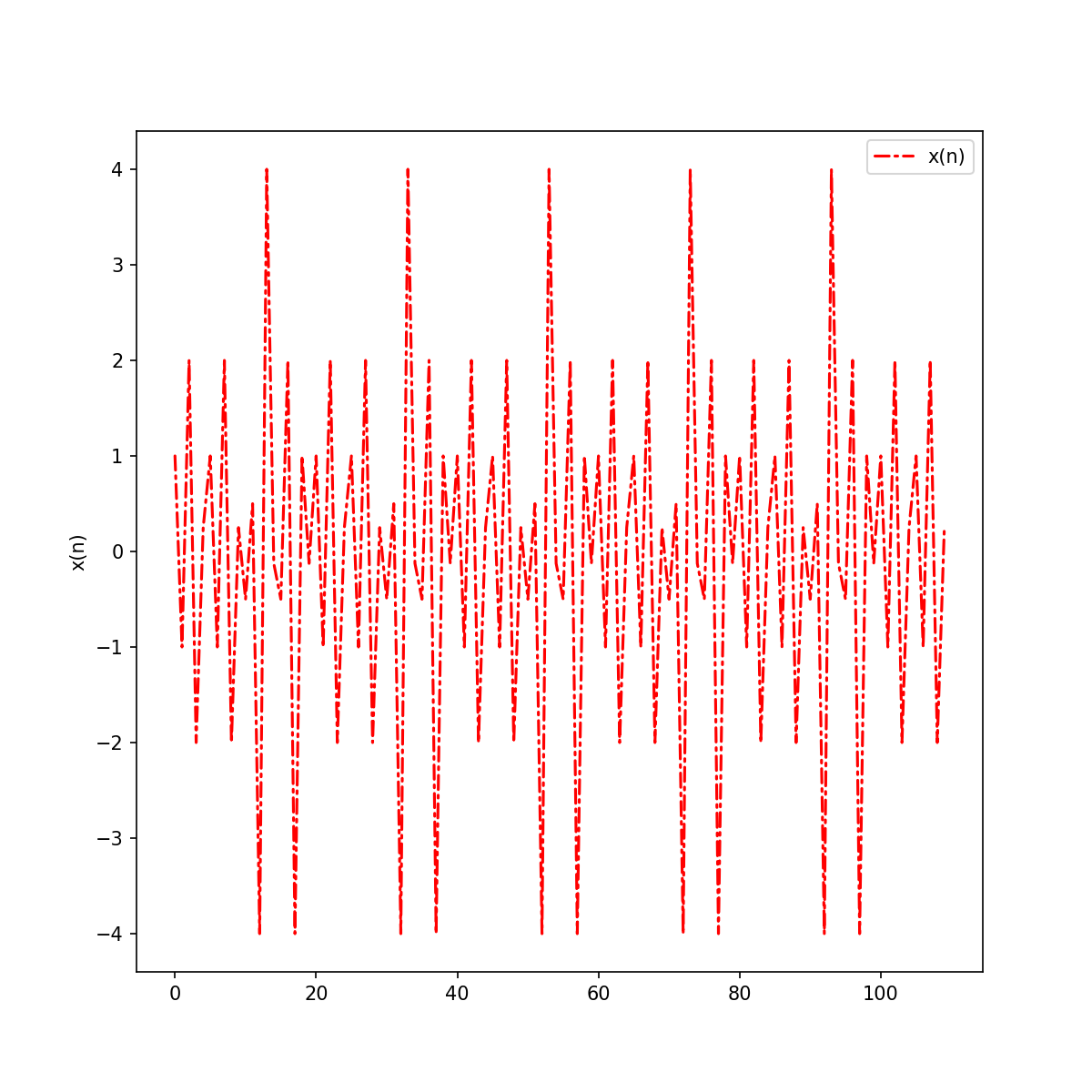}
	\caption{Graph of	{\scriptsize $x_{n+10}=\frac{x_{n}}{-1-x_{n}x_{n+2}x_{n+4}x_{n+6}x_{n+8}}$.}\label{M31_0}}
\end{figure}
In Figure \ref{M31_0}, we use the initial conditions { $x[0] = 1; x[1] = -1; x[2] = 2; x[3] = -2; x[4] = 1/4; 
	x[5] = 1; x[6] = -1; x[7] = 2; x[8] = -2; x[9] = 1/4$}. From the graph, we can verify that the solution is indeed periodic with period $20$.
\begin{Theorem}\label{1}
	Let $x_n$ be a solution of 
	{\footnotesize 
		\begin{equation}\label{eq1cstAnot1}
			x_{n+10}=\frac{x_{n}}{A+Bx_{n}x_{n+2}x_{n+4}x_{n+6}x_{n+8}}
	\end{equation}}
	with initial conditions $x_i, \; i=0,\dots, 9$; and $A\neq 1$ and $B$ some non-zero constants. If the initial conditions   satisfy
	\begin{align}
		&x_{i}x_{i+2}x_{i+4}x_{i+6}x_{i+8}=\frac{1-A}{B}
	\end{align}
	and
	\begin{align}
		& x_i\neq x_{i+2}, x_i\neq x_{i+5},
	\end{align}
	then we have a $10$-periodic solution.
\end{Theorem}
\begin{Proof}
	Assuming $x_{i}x_{i+2}x_{i+4}x_{i+6}x_{i+8}=\frac{1-A}{B}$  as stated in Theorem \ref{1}, equation \eqref{sol1} reduces to 
	\begin{align}\label{I3Anot130p}
		x_{10n+j}=&x_j\prod_{t=0}^{n-1}\frac{A^{5t+\lfloor \frac{j}{2}\rfloor}+{\left(\frac{1-A^{5t+\lfloor \frac{j}{2}\rfloor}}{1-A}\right)B}{\left(\dfrac{1-A}{B}\right)}}{A^{5t+\lfloor \frac{j}{2}\rfloor+1}+{\left(\frac{1-A^{5t+\lfloor \frac{j}{2}\rfloor+1}}{1-A}\right)B}{\left(\dfrac{1-A}{B}\right)}}\\=&x_j.
	\end{align}
	Imposing the condition $ x_i\neq x_{i+2}, x_i\neq x_{i+5}$, we note that the solution is periodic with period $10$. 
\end{Proof}
\begin{figure}[H]
	\centering
	\includegraphics[scale=0.18]{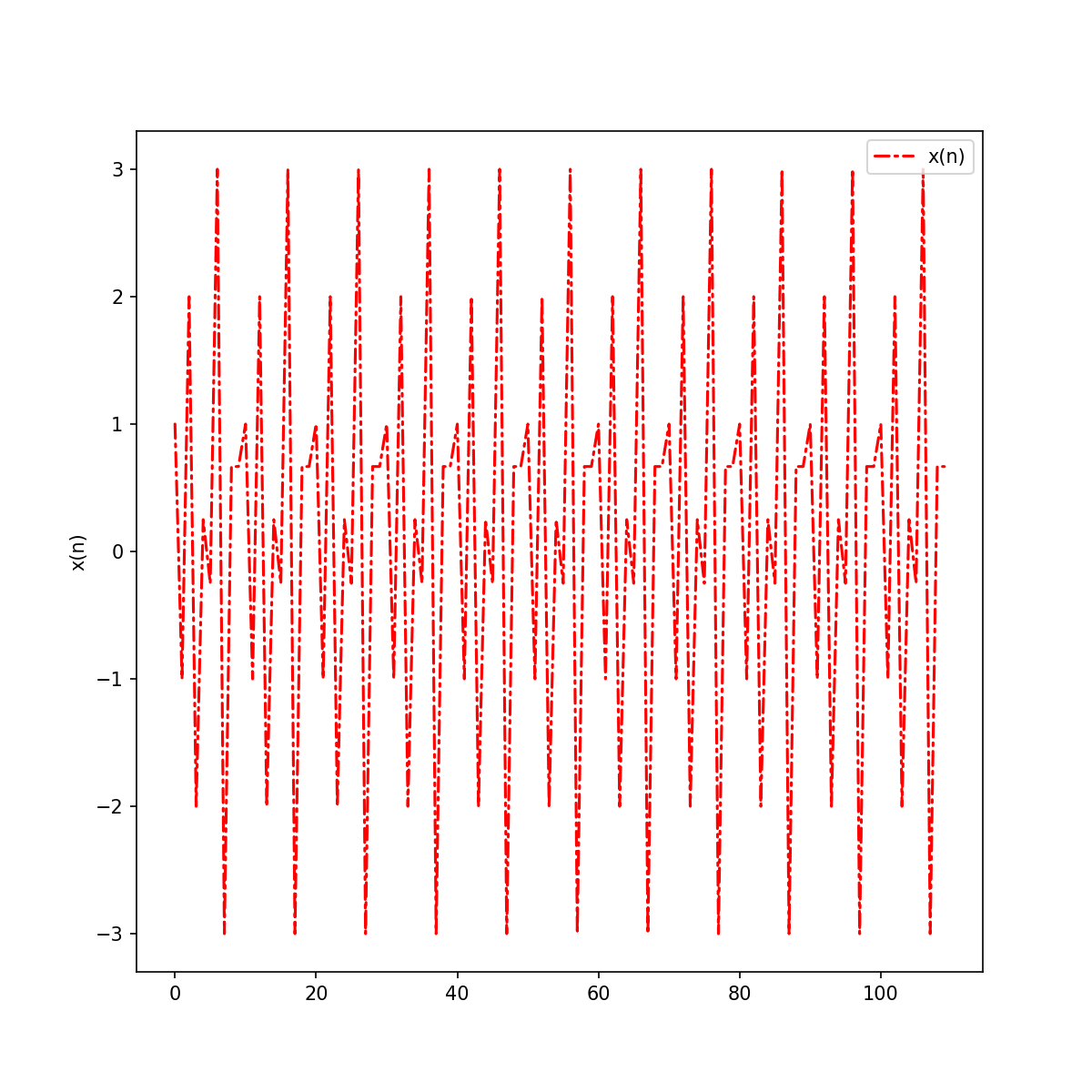}
	\caption{Graph of	{\scriptsize $x_{n+10}=\frac{x_{n}}{2-x_{n}x_{n+2}x_{n+4}x_{n+6}x_{n+8}}$.}\label{M31_1}}
\end{figure}
In Figure \ref{M31_1}, we use the initial conditions { $x[0] = 1; x[1] = -1; x[2] = 2; x[3] = -2; x[4] = 1/4; 
	x[5] = -1/4; x[6] = 3; x[7] = -3; x[8] = 2/3; x[9] = 2/3$} for $A=2$ and $B=-1$. They satisfy
\begin{align}
	x_{0}x_{2}x_{4}x_{6}u_{8}=	x_{1}x_{3}x_{5}x_{7}x_{9}=\frac{1-A}{B}
	\quad \text{and} \quad
	x_i\neq x_{i+2},x_i\neq x_{i+5}.
\end{align}
From the graph, we can verify that the solution is indeed periodic with period $10$.
\begin{Theorem}\label{2}
	Let $x_n$ be a solution of 
	{\footnotesize 
		\begin{equation}\label{eq1cstAnot1'}
			x_{n+10}=\frac{x_{n}}{A+Bx_{n}x_{n+2}x_{n+4}x_{n+6}x_{n+8}}
	\end{equation}}
	with initial conditions $x_i, \; i=0,\dots, 9$; and $A\neq 1$ and $B$ some non-zero constants. If the initial conditions   satisfy
	\begin{align}
		&x_{0}x_{1}x_{2}x_{3}x_{4}=\frac{1-A}{B}
	\end{align}
	and
	\begin{align}
		& x_i\neq x_{i+2}, x_i= x_{i+5},
	\end{align}
	then we have a periodic solution with period $5$.
\end{Theorem}
\begin{Proof}
	The proof is similar to the proof of Theorem \ref{1} and is thus omitted.
\end{Proof}
\begin{figure}[H]
	\centering
	\includegraphics[scale=0.18]{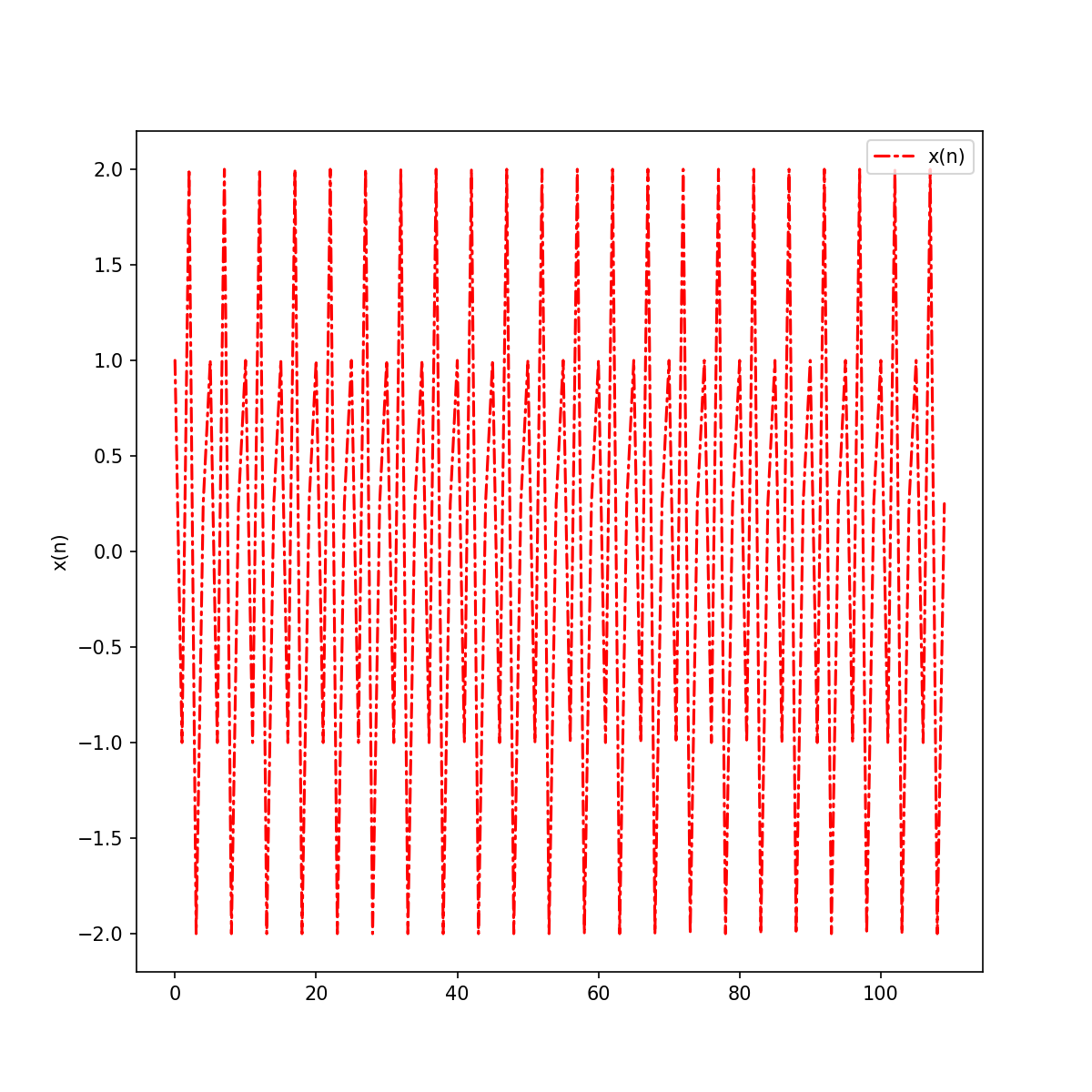}
	\caption{Graph of	{\scriptsize $x_{n+10}=\frac{x_{n}}{2-x_{n}x_{n+2}x_{n+4}x_{n+6}x_{n+8}}$ }.\label{M31_2}}
\end{figure}
In Figure \ref{M31_2}, we use the initial conditions { $x[0] = 1; x[1] = -1; x[2] = 2; x[3] = -2; x[4] = 1/4;
	x[5] = 1; x[6] = -1; x[7] = 2; x[8] = -2; x[9] = 1/4$} for $A=2$ and $B=-1$. They satisfy 
\begin{align}
	x_{0}x_{2}x_{4}x_{6}x_{8}=\frac{1-A}{B}
	\quad 
	\text{and} \quad 
	x_i\neq x_{i+2}, x_i= x_{i+5}.
\end{align} 
From the graph, we can verify that the solution is indeed periodic with period $5$.
\begin{Theorem}\label{3}
	Let $x_n$ be a solution of 
	{\footnotesize 
		\begin{equation}\label{eq1cstAnot1''}
			x_{n+10}=\frac{x_{n}}{A+Bx_{n}x_{n+2}x_{n+4}x_{n+6}x_{n+8}}
	\end{equation}}
	with initial conditions $x_i, \; i=0,\dots, 9$; and $A\neq 1$ and $B$ some non-zero constants. If the initial conditions   satisfy
	\begin{align}
		&x_{i}^5=\frac{1-A}{B}
	\end{align}
	then we have a 1-periodic or constant solution.
\end{Theorem}
\begin{Proof}
	The proof is similar to the proof of Theorem \ref{1} and is thus omitted.
\end{Proof}
\begin{figure}[H]
	\centering
	\includegraphics[scale=0.18]{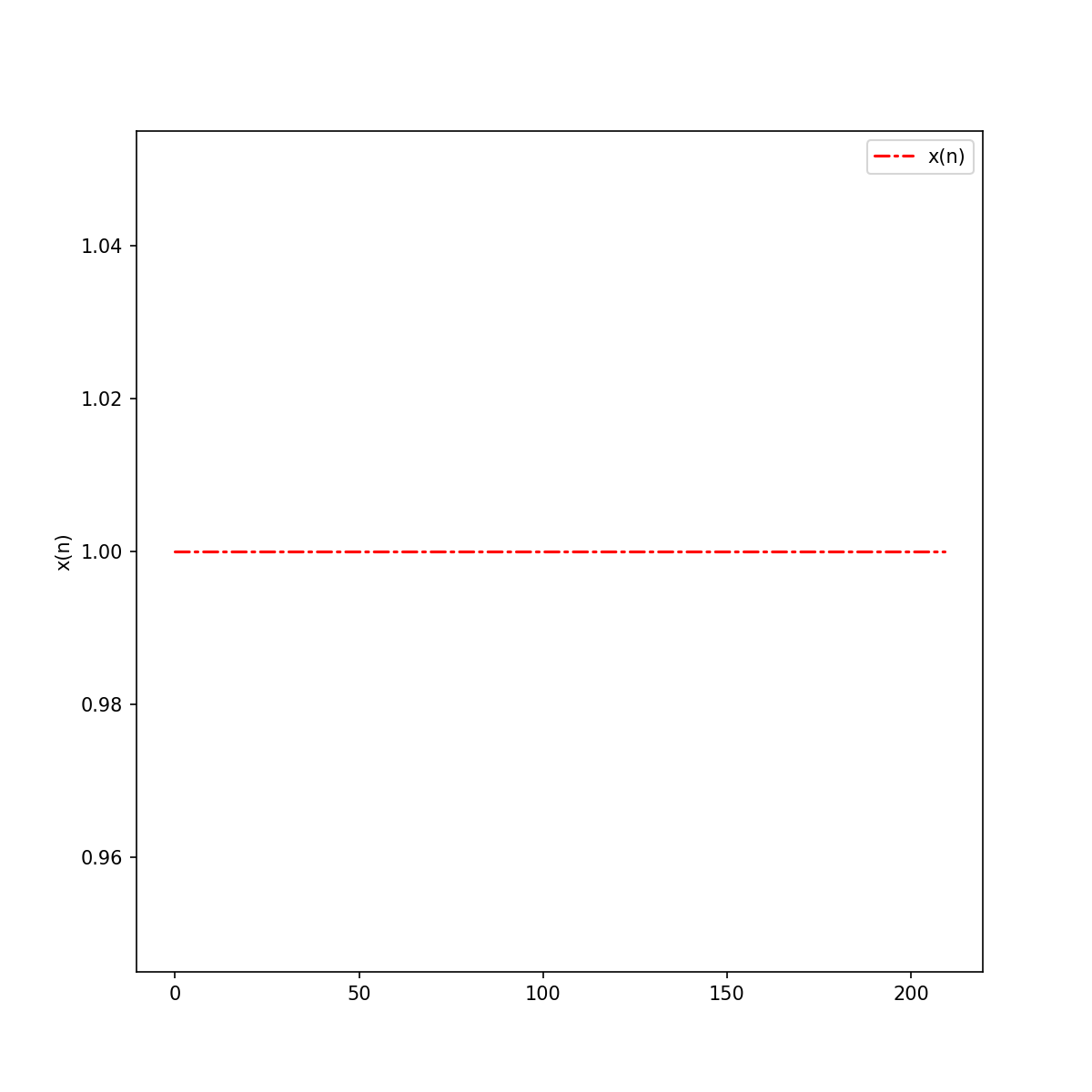}
	\caption{Graph of	{\scriptsize $x_{n+10}=\frac{x_{n}}{2-x_{n}x_{n+2}x_{n+4}x_{n+6}x_{n+8}}$.}\label{M31_3}}
\end{figure}
\begin{figure}[H]
	\centering
		\includegraphics[scale=0.18]{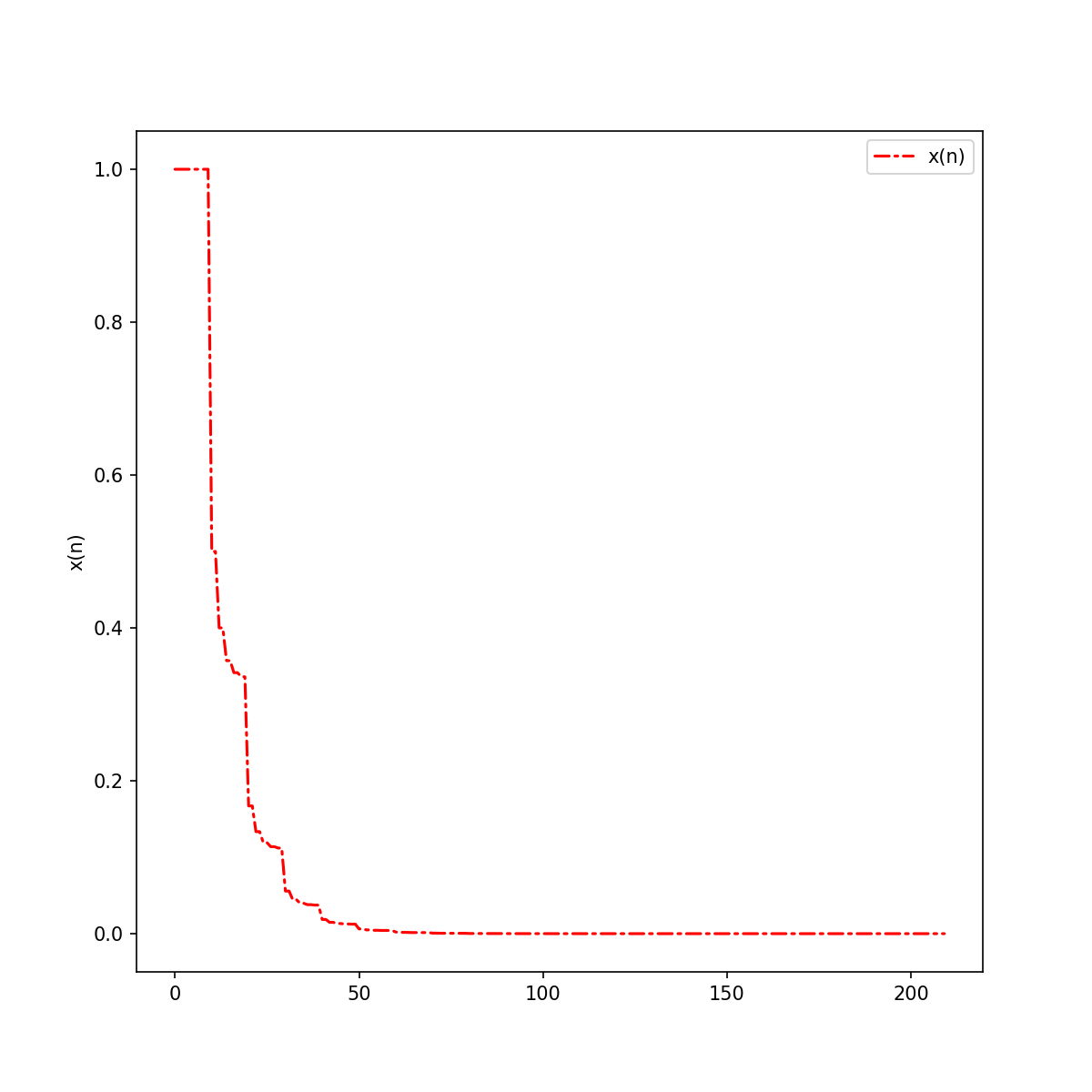}
	\caption{Graph of	{\scriptsize $x_{n+10}=\frac{x_{n}}{3-x_{n}x_{n+2}x_{n+4}x_{n+6}x_{n+8}}$.}\label{M31_4}}
\end{figure}
In Figure \ref{M31_3}, we use the initial conditions  $x[i]=1$ for $A=2$ and $B=-1$. They satisfy
\begin{align}
	x_{i}^5=\frac{1-A}{B}.
\end{align}
From Figure \ref{M31_3}, we can verify that the solution is indeed periodic with period $1$. Despite the fact that all the initial conditions are the same in Figure \ref{M31_4}, where for $A=3$ and $B=-1$, we do not have a $1-$periodic solution because $x_i\neq \frac{1-A}{B}$.
\begin{Theorem}
	Given the equation {\footnotesize
		\begin{equation}\label{eq1cstA1sta}
			x_{n+10}=\frac{x_{n}}{1+Bx_{n}x_{n+2}x_{n+4}x_{n+6}x_{n+8}},
	\end{equation}}
	the only equilibrium point $\bar{x}=0$ is non-hyperbolic.
\end{Theorem}
\begin{Proof}
	The condition $x=x/(1+Bx^5)$ yields the fixed point $x=0$ and the characteristic equation of  \eqref{eq1cstA1sta} near $0$ is $\lambda^{10}-1=0$ whose roots are the tenth roots of unity. The moduli of the roots are all equal to 1. Thus, $x=0$ is  non-hyperbolic.
\end{Proof}
\begin{Theorem}
	The fixed point $\bar{x}=0$ of \eqref{eq1cstAnot1} is (locally) asymptotically stable for $|A|>1$. Furthermore, the non-zero fixed points of \eqref{eq1cstAnot1} are non-hyperbolic for all $A\neq 1$.
\end{Theorem}
\begin{Proof}
	Using the fixed point condition on \eqref{eq1cstAnot1} leads to $\bar{x} (A+B\bar{x}^5-1)=0$. For the first part of the proof, observe that the characteristic equation of \eqref{eq1cstAnot1} near $\bar{x}=0$ is given by {\footnotesize  $\lambda ^{10}-\frac{1}{A}=0$. } Hence, the roots $\lambda_i$ satisfy $|\lambda_i|<1$ for $|A|>1$. Thus, $x=0$ is locally asymptotically stable when $|A|>1$.
	 For the second part  of the proof, we find the non-zero fixed points from the equation $A+B{x}^5-1=0$. In this case, the characteristic equation of \eqref{eq1cstAnot1} near a non-zero equilibrium point  is given by {\footnotesize
		\begin{align}\label{charaanot1}
			0=&\lambda ^{10}-(A-1)\lambda ^{8}-(A-1) \lambda ^{6}- (A-1)\lambda^{4}-{(A-1)} \lambda ^{2}-{A}\\
			&=\frac{1-\lambda^{10}}{1-\lambda^2}( \lambda ^{2}-{A}), \quad \text{with}\quad  \lambda^2\neq1.
	\end{align}}
	There is a root of the characteristic equation (for example, $\lambda= e^{i2\pi/5}$) whose modulus is equal to $1$.
	\end{Proof}
\section{Conclusion}
In this paper, we have confirmed and extended the findings of Sando and Elsayed\cite{SE}. In fact, it has been proved that the findings of this study are consistent with Lama's work through Lie analysis. We have also studied periodicity and stability of the solutions of more general difference equations than the ones considered by Sando and Elsayed \cite{SE}.
\section*{Declarations}

\begin{itemize}
	\item Funding: Not applicable.
	\item Conflict of interest: The authors declare that they have no conflicts of interest.
	\item Ethics approval and consent to participate: Not Applicable.
	\item Consent for publication: Not applicable.
	\item Data availability: Not applicable.
	\item Materials availability: Not applicable.
	\item Code availability: Not applicable.
	\item Author contribution: The authors contributed equally to this work.
\end{itemize}

\end{document}